\documentclass[review]{elsarticle}
\usepackage{lscape}
\usepackage{lineno,hyperref}
\modulolinenumbers[5]
\usepackage{amsmath,amssymb,amsfonts,latexsym}
\DeclareMathOperator{\grad}{grad}

\journal{Journal of \LaTeX\ Templates}

\begin{document}

\begin{frontmatter}

\title{On Time-dependent Hamiltonian Realizations of Planar  and Nonplanar Systems}

\author{O\u{g}ul Esen}
\cortext[mycorrespondingauthor]{Corresponding author}
\address{Department of Mathematics, Gebze Technical University, 41400, Gebze, Kocaeli, Turkey}
\ead{oesen@gtu.edu.tr}

\author{Partha Guha}
\address{S.N. Bose National Centre for Basic Sciences, 
JD Block, Sector III, Salt Lake, Kolkata - 700098, India}
\address{IHES, Le Bois-Marie 35 rue de Chartres 91440 Bures-sur-Yvette France.}
\ead{partha@bose.res.in}

\begin{abstract}
In this paper, we elucidate the key role played by the cosymplectic geometry in the theory of time dependent 
Hamiltonian systems in $2D$. We generalize the cosymplectic structures to time-dependent Nambu-Poisson Hamiltonian
systems and corresponding Jacobi's last multiplier for $3D$ systems.  We illustrate our constructions with various examples.
\end{abstract}


\end{frontmatter}


\tableofcontents

\section{Introduction}

The method of Jacobi's Last Multiplier (JLM) is a geometric way to determine  (possibly non-canonical) Hamiltonian realization of a dynamical system. Here is an incomplete list  \cite{CGK,CGG,Gu13,GCG,Nucci,NuLe04,NuLe08,PoOz17,Leach,Ozie} of the studies on JLM. It is well known that, the method of JLM may result with time-dependent Hamiltonian functions even for autonomous systems. The goal of this work is to investigate this phenomenon for the cases of $2D$ and $3D$ systems. Accordingly, we have organized the main body of the paper into two sections.

In the forthcoming section, we shall focus on $2D$ systems. To investigate the time-dependent cases, we shall address the cosymplectic geometry \cite{Cap,CdLM1,CdLM2,LeSa16,Liber}. As an example, we shall present non-canonical or/and cosymplectic Hamiltonian realizations of a host-parasite model in $2D$. In this example, we shall also employ the theory of conformal Hamiltonian dynamics \cite{MP} to deal with affine terms. The Hamiltonian analysis and JLM diagnosis of some other population growth models in
biological systems will be provided in a list located at the end of this section. 

The first goal of the third section is to elaborate the underlying differential geometric frameworks both of the autonomous and the non-autonomous Nambu-Hamiltonian systems in $3D$. The novelty of this part lying in a proposal of the covariant representation of a time-dependent Nambu-Hamiltonian dynamics (see \ref{NH-3D}). This covariant representation will also lead us for JLM diagnosis of $3D$ systems. To deal with affine terms, the notion of a conformal Nambu-Hamiltonian system will be introduced (see \ref{cNH-Eq}). We shall illustrate these geometries on two chaotic systems, namely L\"{u} and Qi systems.

\section{Hamiltonian systems in Two Dimensions}

\subsection{Canonical Hamiltonian systems}
A manifold $\mathcal{Q}$ is called a symplectic manifold if it is equipped with a closed non-degenerate two-form $\Omega$. 
For a given Hamiltonian (energy) function $H$, the non-degeneracy of the symplectic two-form $\Omega$ manifests existence and uniqueness of the Hamiltonian vector field $X_H$ satisfying the Hamilton's equations
\begin{equation}\label{HamEq}
  i_{X_H}\Omega=dH,
\end{equation}
where $i$ denotes the interior derivative, \cite{dLR,LiMa12,MaRa99}. The dynamics is governed by the Hamiltonian vector field $X_H$.

In Darboux' coordinates $(q^i,p_j)$ on $\mathcal{Q}$, the symplectic two-from can be written as $\Omega=dq^i\wedge p_i$. In this case, the Hamiltonian vector field turns out to be
\begin{equation}\label{HamVecLoc}
X_{H}=\frac{\partial H}{\partial p_i}\frac{\partial }{\partial q^i}-\frac{\partial H}{\partial q^i}\frac{\partial}{\partial p_i}.
\end{equation}
The dynamics is a system of ODEs
\begin{equation}\label{HamEqLoc}
\dot{q}^i=\frac{\partial H}{\partial p_i},\qquad \dot{p}_i=-\frac{\partial H}{\partial q^i}.
\end{equation}

By taking the exterior derivative of (\ref{HamEq}), we arrive at 
\begin{equation}\label{HamPre}
  \mathcal{L}_{X_H}\Omega=0,
\end{equation}
where $\mathcal{L}_X=di_X+i_Xd$ denotes the Lie derivative with respect to $X_H$ \cite{dLR}.
Integration of (\ref{HamPre}) results with the preservation of the symplectic two-from under the flow $\phi_t$ of $X_H$ that is
\begin{equation}\label{symplecticflow}
  \phi_t^*\Omega=\Omega.
\end{equation}
We can argue, more generally, that a vector field is Hamiltonian if and only if its flow preserves the symplectic two-form. Note that, in $2D$, the symplectic two-form can be considered as the area form, this says that a vector field is Hamiltonian if and only if it is divergence free.

\subsection{Conformal Hamiltonian systems} \label{chs}

Let $(\mathcal{Q},\Omega)$ be a symplectic manifold. A vector field $\Gamma^a$ on $\mathcal{Q}$ is called a conformal vector field, if it preserves the symplectic two-form up to some scaling factor $a\in\mathbb{R}$ that is
\begin{equation}\label{conformalvf}
 \mathcal{L}_{\Gamma^a}\Omega =a\Omega,
\end{equation}
see  \cite{Marle}. It is immediate to observe that, a conformal vector field reduces to a Hamiltonian vector field if the scaling factor is zero. The flow $\phi_{t}$ of a conformal vector field preserves the symplectic two-form by the conformal factor $e^{at}$, that is
\begin{equation*}
\phi_{t}^{\ast}\Omega = e^{at}\Omega. 
\end{equation*}

A symplectic manifold $(\mathcal{Q},\Omega)$ admits a conformal vector field with a non-zero scaling factor if and only if the symplectic two-form is exterior derivative of a one-form \cite{MP}. Assume such a case, and take $\Omega = -d\theta$ for some one-form $\theta$. For a given Hamiltonian function $H$, vector field $X_{H}^{a}$  defined through
\begin{equation} \label{CHVF}
i_{X_{H}^{a}}\Omega = dH - a\theta
\end{equation}
is conformal. Locally, the inverse of this assertion is also true. That is, if a vector field is conformal, then there exist a function $H$ and the conformal vector field can be written in the form of $X_{H}^{a}$ satisfying \eqref{CHVF}. More general, the set of conformal vector fields on $\mathcal{Q}$ is
given by $\{X_H + aZ\}$ where $Z$ being the Liouville vector field satisfying $i_Z\Omega = -\theta$.

In Darboux' coordinates $(q^i,p_j)$ on a symplectic manifold $(\mathcal{Q},\Omega)$, the Liouville one-form is $\theta=p_i dq^i$, and the Liouville vector field is $Z=p_i\partial/\partial p_i$. In this local picture picture, a conformal Hamiltonian vector field is computed to be 
\begin{equation}\label{ConfHamVecLoc}
X_H=\frac{\partial H}{\partial p_i}\frac{\partial }{\partial q^i}+(ap_i-\frac{\partial H}{\partial q^i})\frac{\partial }{\partial p_i}.
\end{equation}

\subsection{Cosymplectic manifolds}

A cosymplectic manifold is a triple $({\mathcal{M}},\eta, \Omega)$
consisting of a smooth $(2n+1)-$dimensional manifold ${\mathcal{M}}$ with a closed $1$%
-form $\eta$ and a closed $2$-form $\Omega$ such that $\eta \wedge \Omega^n $ is a non-vanishing volume form on ${\mathcal{M}}$ \cite{Cap,CdLM1,CdLM2}. 

On ${\mathcal{M}}$, there is a distinguished vector field $\xi$, called as Reeb field, determined by the conditions 
\begin{equation} \label{xi}
i_\xi\eta = 1, \qquad i_{\xi}\Omega = 0.
\end{equation} 
There exists an isomorphism $\chi$ from the space ${\mathfrak{X}}({\mathcal{M}})$ of vector fields to the space $\Lambda^1({\mathcal{M}})$ of one-forms that is
\begin{equation}  \label{cosym}
\chi:\mathfrak{X}({\mathcal{M}}) \mapsto \Lambda^1({\mathcal{M}}): X \mapsto \chi(X)= i_X\Omega + \eta(X)\eta.
\end{equation} 
In terms this isomorphism,  the Reeb vector field $\xi$ can be defined by $\chi^{-1}(\eta)$ as well. Accordingly, in the literature, the one-form $\eta$ is also called the Reeb covector field \cite{Cap}.  

The isomorphism $\chi$ presented in (\ref{cosym}) leads us to define a vector field $$\grad H:=\chi^{-1}(dH)$$ on  ${\mathcal{M}}$ for a given real valued function $H$ \cite{CaLeLa92}. This gradient vector field satisfies the following relations 
$$i_{\grad H}\eta =\xi(H) , \qquad i_{\grad H} \Omega = dH-\xi(H)\eta,$$
where $\xi$ is the Reeb vector field.
In Darboux' coordinates $(q^i,p_i,t)$ on ${\mathcal{M}}$, the differential forms are $\Omega=dq^i\wedge dp_i$ and $\eta=dt$. In this local picture, the Reeb vector field $\xi$ and the gradient vector field $\grad H$ take the following forms
\begin{equation}
  \xi=\frac{\partial}{\partial t},\qquad \grad H=\frac{\partial H}{\partial t}\frac{\partial}{\partial t}+\frac{\partial H}{\partial p_i}\frac{\partial }{\partial q^i}-\frac{\partial H}{\partial q^i}\frac{\partial }{\partial p_i}.
\end{equation}

For a Hamiltonian function $H$ on ${\mathcal{M}}$ the associated Hamiltonian vector field $X_H$ is defined through by the equations
\begin{equation} \label{cosymHamVF}
i_{X_H}\eta=0, \qquad i_{X_H}\Omega = dH - \xi(H)\eta.
\end{equation}
The Hamiltonian vector field $X_H$ can also be defined by the preimage of the one-form $dH - \xi(H)\eta$ under the isomorphism $\chi$ as well.
By adding the Reeb field $\xi$ to the Hamiltonian vector field $X_H$, we define the evolution vector field $E_H$. In a compact notation, we define $E_H$ by the following identities 
 \begin{equation} \label{cosymEVF}
i_{E_H}\eta=1, \qquad i_{E_H}\Omega = dH - \xi(H)\eta.
\end{equation}
In Darboux' coordinates $(q^i,p_i,t)$ on ${\mathcal{M}}$, the Hamiltonian vector field and the evolution vector field turn out to be
\begin{equation}\label{local--}
X_H=\frac{\partial H}{\partial p_i}\frac{\partial }{\partial q^i}-\frac{\partial H}{\partial q^i}\frac{\partial }{\partial p_i}, \qquad  E_H=\frac{\partial}{\partial t}+\frac{\partial H}{\partial p_i}\frac{\partial }{\partial q^i}-\frac{\partial H}{\partial q^i}\frac{\partial }{\partial p_i}. 
\end{equation}

\subsection{The method of Jacobi's last multiplier for 2D systems} \label {JLM2D}

Let $\mathcal{Q}$ be a two dimensional manifold with coordinates $(x,y)$. Consider
the following (possibly non-autonomous) system of differential equations
\begin{equation}  \label{e1}
\dot{x}=f(x,y,t),\qquad \dot{y}=g(x,y,t),
\end{equation}
where $f$  and $g$ being real-valued smooth functions. In this section, our main concern is to concentrate on a possible Hamiltonian realization of the system (\ref{e1}), see \cite{Cast} for a similar task. To this end, we first recast the system (\ref{e1}) as the vanishing of the differential one-forms
\begin{equation*}
\alpha^{(1)}:=dx-f(x,y,t)dt,\qquad \alpha^{(2)}:=dy-g(x,y,t)dt,
\end{equation*}
on the extended phase space $\mathcal{M}=\mathcal{Q} \times \mathbb{R}$. We take the exterior product of these one-forms in order to associate the system (\ref{e1}) with the following two-form
\begin{equation} \label{alpha2}
\alpha^{(1)}\wedge\alpha^{(2)}=dx\wedge dy +(fdy-gdx)\wedge dt.
\end{equation}

The Poincar\'e-Cartan one-form for a Hamiltonian $H$ is
given in the canonical coordinates $(q,p,t)$ by
\begin{equation} \label{PC}
\Theta = pdq - Hdt,
\end{equation}
See that, $\Theta$ is consisting of two terms, namely the Liouville one-form $pdq$ and the Hamiltonian term $Hdt$, \cite{dLR}. It is evident that minus of the exterior derivative \begin{equation}\label{e6}
-d\Theta = \beta^{(1)}\wedge \beta^{(2)}=\big(%
dq - {\frac{\partial H}{\partial p}}dt \big)\wedge \big(dp + {\frac{\partial H}{\partial q}}dt \big)
\end{equation}
of $\Theta$ is the exterior product of two one-forms $\beta^{(1)}$ and $\beta^{(2)}$.

If the system (\ref{e1}) has a Hamiltonian realization, then one expects that the two-form $\alpha^{(1)}\wedge\alpha^{(2)}$ presented in (\ref{alpha2}) must be proportional to the closed two-form $\beta^{(1)}\wedge\beta^{(2)}$ in (\ref{e6}). That is, there exists a real valued function $M=M(x,y,t)$, called as the Jacobi's last multiplier, on the extended space $\mathcal{M} = \mathcal{Q} \times \mathbb{R}$ such that
\begin{equation}  \label{e8}
M(dx\wedge dy+(fdy-g dx)\wedge dt)=dq\wedge dp +dH\wedge dt.
\end{equation}
Let us discuss implications of the existence of a Jacobi's last multiplier for a given system. First of all, the two-form on the right hand side of (\ref{e8}) is closed. So that, the two-form on the left hand side must also be closed. This reads the following identity
\begin{equation}  \label{e9}
\left( {\frac{\partial M}{\partial t}}+f{\frac{\partial M%
}{\partial x}} +g{\frac{\partial M}{\partial y}}\right)+M\left({\frac{%
\partial f}{\partial x}}+ {\frac{\partial g}{\partial y}}\right)=0.
\end{equation}
Note that, the second term is nothing but the divergence of the system (\ref{e1}). This is the integrability condition, and can be used to determine $M$. An alternative form of the equation (\ref{e9}) is possible in terms of the total time derivative of $M$ denoted to be
\begin{equation} \label{JLM}
 {\frac{dM}{dt}} +M \hbox{div}X= 0,
\end{equation}
where $X=X(x,y,t)$ is the two dimensional vector field generating the system (\ref{e1}), and $\hbox{div}X$ represents the divergence of $X$ with respect to the space variables $(x,y)$. Here, the total derivative $d/dt$ is the sum $\partial/\partial t + X$.

The form of the defining equation (\ref{JLM}) enables us to arrive at the following observation about the multiplier. If the system $\dot{x}=X$ is divergence free, it is canonically Hamiltonian by choosing $M$ as a constant function, say $M=1$. If not, the Jacobi's last multiplier $M$ must be non-constant. That is, if the generating vector field $X$ is not divergence free then $M$ can not be a constant function. Assume particularly that, $M$ depends only on the space variables, that is $M=M(x,y)$, then the defining equation (\ref{JLM}) takes the particular form 
\begin{equation}\label{JLM-2}
\hbox{div} (M X)= 0.
\end{equation}
In this case, the search of a Jacobi's last multiplier for a given system $\dot{x}=X$ is a search for a function $M$ which makes $MX$ divergence free. 

Secondly, the existence of a Jacobi's last multiplier for a system leads to the determination of a Hamiltonian realization of the system. To see this, assume the existence of a time-independent Jacobi's last multiplier $M$. In this case, we have
\begin{equation}  \label{e10}
M(f dy-gdx)=dH+\mbox{terms proportional to }\; dt.
\end{equation}
In general, multiplicative inverse of the multiplier is observed as one of the coefficients of the symplectic structure \cite{Gu13}. 
Consider the following non-canonical Poisson structure 
\begin{equation}\label{npr}
  \{x, x\} = \{y, y\} = 0, \qquad \{x,y\}= \frac{1}{M} .
\end{equation}
on the two-dimensional phase space $\mathcal{Q}$. If $M$ is non vanishing, then this Poisson structure becomes non-degenerate. In this case, we can write the symplectic two-form as
\begin{equation}\label{npr-}
  \Omega = \frac{1}{M} dx\wedge dy.
\end{equation}
Here, $ {1}/{M} $ plays the role of a conformal factor. 
Using this symplectic structure, we write the Hamilton's equations as
\begin{equation}\label{nsHe}
  \dot{x} = \frac{1}{M} {\frac{\partial H}{\partial y}}, \qquad \dot{y} = - \frac{1}{M}  {%
\frac{\partial H}{\partial x}}.
\end{equation}

For time-dependent cases, a generalization of this geometry is possible by modifying (\ref{e10}) properly. To this end, introduce two auxiliary functions $\psi$
and $\phi$  such that
\begin{equation}  \label{e11}
M((f-\psi) dy-(g-\phi)dx)=dH+\mbox{terms proportional to}\; dt.
\end{equation}
This essentially removes
the explicit time-dependent terms and allows for the construction of a
Hamiltonian for the remaining autonomous part \cite{Cast}. However, the time dependence
is not altogether lost, instead it is manifested in the the data of the new
coordinates. Now consider the following reduced system after the subtraction of the auxiliary parts
\begin{equation}\label{sys2}
\begin{cases}
\dot{x}=f(x,y,t)-\psi(x,y,t),\\ 
\dot{y}=g(x,y,t)-\phi(x,y,t).
\end{cases}
\end{equation}
This reduced system is Hamiltonian if we introduce the non-canonical Poisson relations in (\ref{npr}). Now, the divergence free criteria (\ref{JLM-2}) takes the form of
\begin{equation}  \label{e12}
{\frac{\partial}{\partial x}}\left(M(f-\psi)\right) +{\frac{%
\partial}{\partial y}}\left(M(g-\phi)\right)=0.
\end{equation}
In this generalized framework, the initial defining equation (\ref{e8}) of the Jacobi's last multiplier should be modified as
\begin{eqnarray} \label{Mg}
M(dx\wedge dy +(f dy-g dx)\wedge dt)&=& M (dx-\psi
dt)\wedge(dy-\phi dt)  \\ && +M((f-\psi)dy\wedge dt-(g-\phi)dx\wedge dt) \notag  \\
&=&M(dx-\psi dt)\wedge
(dy-\phi dt)+dH\wedge dt, \notag 
\end{eqnarray}
where we have employed the equation (\ref{e11}) while arriving at the last line of the calculation. 

On the other hand, a direct comparison of the right hand side of (\ref{e8}) and the calculation (\ref{Mg}) results with the following equation
\begin{equation}  \label{e14}
M(dx-\psi dt)\wedge (dy-\phi dt)=dq\wedge dp.
\end{equation}
This enables us to find a (possibly time-dependent) symplectic transformation determining the canonical coordinates of the physical system. Thus the problem of recasting (\ref{e1}) into the form
of Hamilton's equations reduces to the determination of auxiliary
functions $\phi$ and $\psi$ such that $H$ is identified from (\ref{e11}). Note that, the equation (\ref{e14}) brings additional integrability conditions on the auxiliary variables and the Jacobi's last multiplier if the determinations of the canonical coordinates are insisted. This may results with time-dependent Hamiltonian functions even for autonomous systems.

Here is the particular form of the systems studied in this present study. Consider an autonomous system under the existence of the Jacobi's last multiplier $M=M(x,y,t)$ depending explicitly on the time variable, and two time-independent auxiliary functions $\psi$ and $\phi$ satisfying (\ref{e14}). The Hamiltonian realization of the induced system 
\begin{equation}  \label{e1-}
\begin{cases}
\dot{x}=M(x,y,t)(f(x,y)-\psi(x,y)), \\ \dot{y}=M(x,y,t)(g(x,y)-\phi(x,y))
\end{cases}
\end{equation}
is possible in the local form (\ref{HamEqLoc}). But, such a system fails to satisfy the Hamilton's equations in the global form (\ref{HamEq}) since the time dependency of the system manifests the time dependency of the Hamiltonian function. To solve this, we address the (cosymplectic) Hamilton's equations presented in (\ref{cosymHamVF}).

Consider a two dimensional phase space $\mathcal{Q}$ with local coordinates $(x,y)$ with the symplectic two-form $\omega=dx\wedge dy$. Consider the canonical inclusion of $\mathcal{Q}$ into the three dimensional product manifold $\mathcal{M}=\mathcal{Q}\times \mathbb{R}$. There exists a projection $\pi$ from $\mathcal{M}$ to $\mathcal{Q}$ which makes them a fiber bundle. This structure enables us to pull the symplectic two-form $\Omega$ to the total space $\mathcal{M}$. We choose the Lee one-form as $\eta=dt$ and the Reeb vector field as $\xi=\partial / \partial t$, then $(\mathcal{M},dt,\pi^*(\Omega))$ becomes a cosymplectic manifold. In this construction, we consider the evolution vector field 
\begin{equation*}
  E_H (x,y,t) = \frac{\partial}{\partial t}+X_H={
\frac{\partial}{\partial t}} + {\frac{\partial H(x,y,t)}{\partial y}}{\frac{\partial%
}{\partial x}} - {\frac{\partial H(x,y,t)}{\partial x}}{\frac{\partial}{\partial y}},
\end{equation*}
hence the dynamics governed by the Hamiltonian vector field $X_H$ can be written as in the coordinate free form (\ref{cosymHamVF}).

\subsection{Example: A host-parasite model} \label{hpm}

We consider the following host-parasite model described by the system of differential equations
\begin{align}  \label{hp.1}
\begin{cases}
\dot{x}& =a x- b yx \\
\dot{y}&= c y- \delta y^2/x
\end{cases}
\end{align}
where $a,b,c,\delta$ are some real constants. We take the auxiliary functions $\phi=ax$, $\psi=cy$, and the Jacobi's last multiplier
\begin{equation}  \label{hp.2}
M(x,y,t)={\frac{e^{ct}}{xy^2}},
\end{equation}
so that from (\ref{e11}) we have
\begin{equation}  \label{hpexact}
{\frac{e^{ct}}{xy^2}}\left\{(a-by)x-\psi\right\}dy- {\frac{e^{ct}}{xy^2}}%
\left\{y(c-\delta{\frac{y}{x}})-\phi\right\}dx=dH.
\end{equation}
See that, these choices lead us that the left hand side of (\ref{hpexact}) is exact. This allows us 
to obtain the Hamiltonian function
\begin{equation}  \label{hp.4}
H=\left(-b\ln y-{\frac{\delta}{x}}\right)e^{ct}.
\end{equation}

Accordingly, to exhibit the cosymplectic Hamiltonian character of the host-parasite model in terms of $(x,y,t)$ coordinates we recall the definition of the evolutionary vector field $E_H$. Define a two form $\pi^*(\omega) = dx \wedge dy$, a one-form $\eta=dt$, and the Reeb vector field $\xi=\partial / \partial t$. Hence a direct calculation shows that the system (\ref
{hp.1}) can be written by 
\begin{equation*}
i_{E_H}\pi^*(\omega) = {\frac{\delta}{x^2}}e^{ct}dx - {\frac{b}{y}} e^{ct} dy =
dH - cHdt \equiv dH - \xi(H)dt
\end{equation*}
which obeys the definition (\ref{cosymHamVF}). 

In order to find the canonical coordinates related with the system (\ref{hp.1}), we recall the definition in (\ref{e14}). So that we arrive at the following identities
\begin{equation*}
dq\wedge dp={\frac{e^{ct}}{xy^2}}(dx-ax dt)\wedge (dy-cy dt)=d(at-\ln x)\wedge d\left({\frac{e^{ct}}{y}}\right).
\end{equation*}
which enable us to determine the canonical coordinates as $q=at-\ln x$ and $p=e^{ct}/y$. In terms of the canonical coordinates, the Hamiltonian function in (\ref{hp.4}) has the appearance
\begin{equation}  \label{hp.6}
H(q, p, t)=-e^{ct}\left[b(ct-\ln p )+\delta e^{q-at}\right],
\end{equation}
whereas the Hamilton's equations become
\begin{equation*}
\dot{q}={\frac{\partial H}{\partial p}}=b{\frac{e^{ct}}{p}},\qquad \dot{p}=-{%
\frac{\partial H}{\partial q}}=\delta e^{q}e^{(c-a)t}.
\end{equation*}

We now present how one can write the Host-Parasite system (\ref{hp.1}) in the form (\ref{CHVF}) of a conformal Hamiltonian system. To this end, we first  set $c = 0$ in (\ref{hp.1}) and obtain the following reduced system
\begin{align}  \label{hp.chr}
\begin{cases}
\dot{x}& =a x- b yx \\
\dot{y}&= - \delta y^2/x
\end{cases}
\end{align}
In this present case, the Jacobi's last multiplier reduces to $M=1/xy^2$ without introducing any auxiliary functions. The Hamiltonian function is computed to be
\begin{equation}\label{HamHP2-}
  H(x,y)=-\frac{a}{y}-b\ln y-{\frac{\delta}{x}}.
\end{equation}
Here, $x$ and $y$ cannot be zero. 
Let us recall the Poisson bracket in (\ref{npr}) and write the corresponding symplectic two-form for this case as follows
\begin{equation}\label{HamHP2--} 
 \Omega=xy^2 dx\wedge dy.
\end{equation}
It is now a matter of direct calculation to show that the reduced system (\ref{hp.chr}) is Hamiltonian
generating by the Hamiltonian function (\ref{HamHP2-}) and the symplectic two-form $\Omega$ in (\ref{HamHP2--}). There remains to add the term $cy$ to the right hand side of the second equation in (\ref{hp.chr}). We will add this term in a pure geometrical way, by introducing the Liouville vector field $Z=y\partial / \partial y$. Then the system (\ref{hp.1}) is in the conformal Hamiltonian form according to the representation
$$i_{(X_H+cZ)}\Omega = dH-c\theta, $$
where $\Omega$ is the symplectic form in (\ref{HamHP2--}) and that $\theta$ is the associated potential one-form $\theta=\frac{1}{3}y^3xdx$ satisfying $\Omega=-d\theta$.

\subsection{Some other examples} \label{table}
\smallskip

We can repeat similar calculations for some other population growth models in
biological systems to obtain their integrating factors, Hamiltonian and
canonical coordinates. In the following table we encode these results. One
must note that the integrating factors of the first three equations are
time-dependent and the next two are time-independent. The first three models
describe cosymplectic structure and their integrating factors are explicitly
time-dependent, and the last two are connected to conformal Hamiltonian
systems. Of course if we kill the time-dependent part for a special value of
the parameters all of them satisfy conformal Hamiltonian structure.

\newpage 

\begin{landscape}

\begin{tabular}{c|c|c|c|c} 
\hline
&  &  &  &  \\
Equations & Models & I.F. & Hamiltonians & Canonical coords \\ \hline\hline
&  &  &  &  \\
$\dot{x} =(a-by)x$ & Host-Parasite & ${\frac{e^{ct}}{xy^2}}$ & $H= (-b\ln y-{%
\frac{\delta}{x}})e^{ct}$ & $q=at-\ln x$ \\
$\dot{y}= (c-\delta {\frac{y}{x}})y$ &  &  &  & $p={\frac{e^{ct}}{y}}$ \\
&  &  &  &  \\
$\dot{x} =x (a\log ({\frac{x}{\kappa}})+by)$ & Gompertz & $e^{-(a+c)t}$ & $H
=e^{-(a+c)t}(by-\delta x)$ & $q=e^{-at}\log ({\frac{x}{\kappa}})$ \\
$\dot{y}=y\left(\delta x+c\log\left({\frac{y}{\nu}}\right)\right)$ &  &  &
& $p=e^{-ct}\log\left({\frac{y}{\nu}}\right)$ \\
&  &  &  &  \\
$\dot{x}= x(\lambda - ax + by),$ & mutualistic & $e^{\alpha t} x^{\beta}
y^{\gamma}$ & $H = -k e^{\alpha t}x^{\beta + 2} y^{\gamma + 1}$ & $q = {%
\frac{x^{\beta+1}e^{-\lambda(\beta + 1)t}}{\beta + 1}}$ \\
$\dot{y} = y(\nu + cx - dy)$ &  &  & $+ k^{\prime} e^{\alpha t} y^{\gamma +
2}x^{\beta + 1}$ & $p = {\frac{y^{\gamma+1}e^{-\nu(\gamma + 1)t}}{\gamma + 1}%
}$ \\
&  &  &  &  \\
$\dot{x}= -x + {\frac{x^2}{y}},$ & Koch & ${\frac{1}{x^2}}$ & $H = \ln y -x $
& $q=t+\ln\sqrt{xy}-{\frac{y}{x}}$ \\
$\dot{y} = -y + x^2$ & -Meinhardt &  &  & $p=t+\ln\sqrt{xy}$ \\
&  &  &  &  \\
$\dot{x}= -\beta xy - \nu x,$ & Kermeck & ${\frac{1}{xy}}$ & $H= -\beta(x+y)$
& $q = \ln y - \nu t$ \\
$\dot{y} = \beta xy - (\nu + \gamma)y$ & McKendrick &  &  & $p = \ln y +
(\gamma + \nu)t$ \\
&  &  &  &  \\ \hline
\end{tabular}

\end{landscape}

\newpage

\section{Nambu-Hamiltonian Systems in Three Dimensions}
\subsection{Poisson manifolds in three dimensions}

Poisson bracket on an $n$-dimensional manifold $\mathcal{P}$ is a binary operation
$\{\bullet,\bullet\}$ on the space $\mathcal{F}(\mathcal{P})$ of real-valued smooth functions satisfying
the Leibniz and the Jacobi identities \cite{LaPi12,LiMa12,wei83}. In this picture, a
system of ODE's is called a Hamiltonian system if it can be written in the form of Hamilton's
equation%
\begin{equation}\label{HamEqn}
\dot{x}=\left\{x,H\right\}%
\end{equation}
for a Hamiltonian function $H$.
A dynamical system is called bi-Hamiltonian if it admits two different Hamiltonian
structures
\begin{equation}
\dot{x}= \{\dot{x},H_2\}_{1}=\{\dot{x},H_1\}_{2},\label{biHam}%
\end{equation}
with the requirement that the Poisson brackets $\{\bullet,\bullet \}_{1}$ and $\{\bullet,\bullet \}_{2}$ be
compatible \cite{Fe94,OLV}.

For a three dimensional Euclidean space, the Hamilton's equation (\ref{HamEqn}) takes
the particular form%
\begin{equation} \label{HamEq3}
\dot{x}=\mathbf{J}\times\nabla H_2,
\end{equation}
\cite{EsGhGu16, Gum1, GuNu93,Raz12,Raz12b}. Here, $\mathbf{J}$ is the Poisson vector corresponding to the skew-symmetric Poisson tensor \cite{MaRa99}. In this case, the Jacobi identity becomes
\begin{equation} \label{jcbv}
\mathbf{J}\cdot(\nabla\times\mathbf{J})=0.
\end{equation}
The general solution of the Jacobi identity (\ref{jcbv}) is given by $(1/M)\nabla H$ for two arbitrary functions $M$ and $H$ \cite{AGZ,HB1,HB2,HB3}. Here, the existence of a Jacobi's last multiplier $M$ is a manifestation of  conformal invariance
of the Jacobi identity \cite{Gor01,Jac1,Jac2,Whi}. This enables us to arrive at the following assertion. A Hamiltonian system in three dimensions is bi-Hamiltonian in the form
\begin{equation} \label{bi-Ham}
\dot{x}=\frac{1}{M}\nabla H_{1}\times\nabla H_{2}=\mathbf{J}_{1}\times\nabla H_{2}=\mathbf{J}_{2}\times\nabla
H_{1}.
\end{equation}
Here, the first Poisson vector $\mathbf{J}_{1}$ is given by $(1/M)\nabla H_{1}$ whereas the second Poisson vector $\mathbf{J}_{2}$ is given by $-(1/M)\nabla H_{2}$. 

\subsection{Three dimensional Nambu-Poisson manifolds}

Let $\mathcal{P}$ be manifold. A Nambu-Poisson bracket of order $3$ is a multilinear ternary operation, denoted by $\{\bullet ,\bullet ,\bullet
\}$, on the space $\mathcal{F}(\mathcal{P})$ of smooth functions. A Nambu-Poisson structure must satisfy both the generalized
Leibniz identity
\begin{equation} \label{GLI}
\left\{ F_{1},F_{2},FH\right\} =\left\{ F_{1},F_{2},F\right\} H+F\left\{
F_{1},F_{2},H\right\}
\end{equation}%
and the fundamental (or Takhtajan) identity
\begin{equation}
\left\{ F_{1},F_{2},\{H_{1},H_{2},H_{3}\}\right\}
=\sum_{k=1}^{3}\{H_{1},...,H_{k-1},\{F_{1},F_{2},H_{k}\},H_{k+1},...,H_{3}\},
\label{FI}
\end{equation}%
for arbitrary functions $F,F_{1},F_{2},H,H_{1},H_{2}$ \cite{Nambu, Ta}. See also \cite{AlGu96,Na98}. 

Assume that $(\mathcal{P},\{\bullet ,\bullet ,\bullet
\})$ be a Nambu-Poisson manifold. For a pair $(H_1,H_2)$ 
of Hamiltonian functions, the associated Nambu-Hamiltonian vector field $X_{H_1,H_2}$ is defined through 
\begin{equation}\label{NambuHamVF}
  X_{H_1,H_2}(F)=\{F,H_1,H_2\}.
\end{equation}
The distribution of the Nambu-Hamiltonian vector fields are in involution and defines a foliation of the manifold $\mathcal{P}$. A
dynamical system is called Nambu-Hamiltonian with a pair $(H_1,H_2)$ of Hamiltonian
functions if it can be recasted as%
\begin{equation}
\dot{x}=\left\{x ,H_{1},H_{2}\right\} .  \label{NHamEqn}
\end{equation}

By fixing one of the Hamiltonian functions in the pair $(H_1,H_2)$, we can write a Nambu-Hamiltonian system in the bi-Hamiltonian form as well
\begin{equation} 
\dot{x}=\left\{ x,H_{1}\right\} ^{H_{2}}=\left\{ x
,H_{2}\right\} ^{H_{1}}  \label{NH-2Ham}
\end{equation}%
Here, the brackets $\{\bullet ,\bullet \}^{H_{2}}$ and $\{\bullet ,\bullet
\}^{H_{1}}$ are compatible Poisson structures defined by
\begin{equation}
\left\{ F,H\right\} ^{{H}_{2}}=\left\{ F,H,H_{2}\right\} \text{, \ \ \ }%
\left\{ F,H\right\} ^{{H}_{1}}=\left\{ F,H_{1},H\right\} ,  \label{Pois}
\end{equation}%
respectively.

Let $N$ be a three-vector field on a manifold $\mathcal{P}$ and define the following multilinear ternary bracket
\begin{equation}\label{N3V}
  \left\{ F_{1},F_{2},F_3\right\}=N(dF_1,dF_2,dF_3),
\end{equation}
where $dF_i$ stands for the de Rham exterior derivatives of the functions $F_i$, for $i=1,2,3$. By definition, this bracket satisfies the skew-symmetry and Leibnitz identity (\ref{GLI}) so that it is a generalized almost Poisson bracket \cite{IbLeMaDi97}. Such a bracket does not necessarily satisfy the Takhtajan identity (\ref{FI}). A three-vector $N$ defines a mapping $\sharp$ from the space $\Lambda^2(Q)$ of two-forms on $\mathcal{P}$ to the space $\mathfrak{X}(\mathcal{P})$ of vector fields on $\mathcal{P}$ by
\begin{equation}\label{sharp}
\langle \sharp(\alpha\wedge\beta),\gamma \rangle = N(\alpha,\beta,\gamma).
\end{equation}
The bracket defined in (\ref{N3V}) is a Nambu-Poisson bracket if a vector field $X_{H_1,H_2}$, defined as \begin{equation}\label{NambuHVF2}
  X_{H_1,H_2}=\sharp(dH_1\wedge dH_2),
\end{equation}
is a derivation on the algebra that is if
 \begin{eqnarray}\label{Int}
 X_{H_1,H_2} \{F_1,F_2,F_3\}&=&\{X_{H_1,H_2} F_1,F_2,F_3\}+\{ F_1,X_{H_1,H_2}F_2,F_3\}\notag \\&&+ ~\{F_1,F_2,X_{H_1,H_2} F_3\},
 \end{eqnarray}
for all $F_1,F_2,F_3$. 
In this case, we call $N$ as a Nambu-Poisson three-vector field. Note that, two definitions of the Nambu-Hamiltonian vector fields in (\ref{NambuHamVF}) and (\ref{NambuHVF2}) are coinciding.

\subsection{Three dimensional Nambu-Hamiltonian systems}
Let $\mathcal{P}$ be a 3 dimensional manifold equipped with a non-vanishing volume manifold $\mu$. Then the following identity 
\begin{equation}\label{NP-Vol} 
  \{F_1,F_2,F_3\}\mu=dF_1\wedge dF_2 \wedge dF_3
\end{equation}
defines a Nambu-Poisson bracket on $\mathcal{P}$ \cite{Ga96,Gu01}. In this case, the distribution defined by the Hamiltonian vector fields is three dimensional. This gives that corresponding  foliation of the distribution consists of a unique leaf \cite{IbLeMaDi97,LeSa17}. In other words, the kernel of $\sharp$ operator defined in (\ref{sharp}) is trivial. 

Assume that the Nambu-Poisson bracket on $\mathcal{P}$ is defined by means of a volume form as described in (\ref{NP-Vol}). In this case, the equation (\ref
{NambuHamVF}) relating a Hamiltonian pair $(H_1,H_2)$ with a Nambu-Hamiltonian vector field $X_{H_1,H_2}$ can be written, in a covariant formulation, as 
\begin{equation} \label{NHcf}
i_{X_{H_1,H_2}}\mu = dH_1 \wedge dH_2.
\end{equation}
We shall call (\ref{NHcf}) as Nambu-Hamilton's equations \cite{Fe92}.
Note that, by taking the exterior derivative of both hand side of (\ref{NHcf}), we arrive at the preservation of the volume form by the Nambu-Hamiltonian vector field, that is
\begin{equation} \label{NH-vf-mu}
\mathcal{L}_{X_{H_1,H_2}}\mu =0.
\end{equation}
Integration of this conservation law gives that the flows of Nambu-Hamiltonian vector fields are volume preserving diffeomorphisms.

Consider a local frame (called as the standard basis) given by a three-tuple $(u,v,w)$ such that the volume form is 
\begin{equation}\label{volume}
 \mu= du \wedge dv \wedge dw.
\end{equation}
In this picture the Nambu-Poisson three-vector (\ref{N3V}) takes the particular form
\begin{equation}\label{Nambu}
N= \frac{\partial}{\partial u} \wedge  \frac{\partial}{\partial v}  \wedge  \frac{\partial}{\partial w}.
\end{equation}
Locally, the Nambu-Hamiltonian vector field $X_{H_1,H_2}$ defined in (\ref{NHcf}) for the pair $(H_1,H_2) $ of Hamiltonian functions can be computed as
\begin{equation} \label{HbhVF}
X_{H_1,H_2} =  \{H_1,H_2\}_{u,v} {\frac{%
\partial}{\partial w}} + \{H_1,H_2\}_{v,w} {\frac{\partial}{\partial u}} + 
\{H_1,H_2\}_{w,u} {\frac{\partial}{\partial v}},
\end{equation}
where the coefficient functions are computed to be, for example,
\begin{equation} \label{pbh}
\{H_1,H_2\}_{a,b} = {\frac{\partial H_1}{\partial a}}{\frac{\partial H_2}{%
\partial b}} - {\frac{\partial H_1}{\partial b}}{\frac{\partial H_2}{%
\partial a}}.
\end{equation}

For the particular case of a three dimensional Euclidean space, the present discussion reduces to the following form. Let $F_1$, $%
F_{2}$ and $F_{3}$ be three real valued functions, and consider the triple product
\begin{equation}
\left\{ F_{1},F_{2},F_{3}\right\} = \nabla F_1 \cdot \nabla F_{2}\times
\nabla F_{3}  \label{NambuPois}
\end{equation}%
of the gradients of these functions. It is evident that the bracket \eqref{NambuPois} is a Nambu-Poisson bracket with corresponding Nambu--Poisson three--vector field in the standard form (\ref{Nambu}). The Nambu-Hamiltonian vector field presented in (\ref{HbhVF}) takes the particular form
\begin{equation*}
X_{H_1,H_2} = \nabla H_{1}\times\nabla H_{2}.
\end{equation*}
It follows that, the Nambu-Hamilton's equations (\ref{NHamEqn}) turn out to be
\begin{equation} \label{NH3D}
\dot{x}=\left\{x ,H_{1},H_{2}\right\}= \nabla H_{1}\times\nabla H_{2}.%
\end{equation}
The bi-Hamiltonian character of this system can easily be observed by employing (\ref{NH-2Ham}).

\subsection{Non-autonomous Three dimensional Nambu-Hamiltonian systems} \label{Na3D}

Now, we establish a geometrical framework of non-autonomous $3$ dimensional Nambu-Hamiltonian systems. To this end, consider a $3$ dimensional volume (hence according to (\ref{NP-Vol}) a Nambu-Possion) manifold $(\mathcal{P},\mu)$. Consider the following product manifold $\mathcal{P}\times \mathbb{R}$ which is endowed with a local coordinate system $(u,v,w,t)$. It is possible to define two projections from the product space to its components. The first projection $\pi_1$ is a surjective mapping from the total space $\mathcal{P}\times \mathbb{R}$ to the real numbers $\mathbb{R}$. Consider a basis $dt$ for the module of differential one-forms on $\mathbb{R}$ and pull this one-form back to $\mathcal{P}\times \mathbb{R}$ by means of $\pi_1$. This results with a horizontal one-form $\eta$ with respect to the fibration $\pi_1$.  Note that, in the local chart $(u,v,w,t)$, this horizontal one form $\eta$ can be written as $dt$ as well. The second projection $\pi_2$ is a surjective mapping from the total space $\mathcal{P}\times \mathbb{R}$ to the volume manifold $\mathcal{P}$. Pull-back of the volume form $\mu$ on $\mathcal{P}$ back to the total space by means of the projection $\pi_2$ is a well-defined constant three-form on $\mathcal{P}\times \mathbb{R}$. For the sake of the clearance of the notation, we shall not distinguish three forms $\mu$ and $\pi_2^*\mu$, and denote both of them by $\mu$. It is easy to observe that the exterior product of $\mu$ and $\eta$ is a non-vanishing top-form on the total space $\mathcal{P}\times \mathbb{R}$. Motivating from the definition of the Reeb vector field presented in (\ref{xi}), we define a vector field $\nu$ on $\mathcal{P}\times \mathbb{R}$ via the following identities
\begin{equation} \label{xiNP}
\qquad i_{\nu}\mu = 0, \qquad i_\nu\eta = 1.
\end{equation}
See that the vector field $\nu$ is a vertical vector field for the projection $\pi_2$. In terms of the local coordinates $(u,v,w,t)$, $\nu$ is written by $\partial / \partial t$. 

We define evolutionary vector field $E_{H_1,H_2}$ associated with a (possibly time-dependent) Hamiltonian function pair $(H_1,H_2)$ by means of the following equalities
\begin{equation}\label{tNH}
  i_{E_{H_1,H_2}}\mu=dH_1\wedge dH_2 - \nu(H_1)\eta\wedge dH_2 - dH_1\wedge \nu(H_2)\eta, \qquad i_{E_{H_1,H_2}}\eta = 1,
\end{equation}
where $\eta$ is the horizontal one-form, and $\nu$ is the vector field in (\ref{xiNP}). Locally the equation (\ref{tNH}) reads the Evolution $E_{H_1,H_2}$ in form
\begin{eqnarray}
E_{H_1,H_2}= {\frac{\partial}{\partial t}} +  \{H_1,H_2\}_{u,v} {\frac{%
\partial}{\partial w}} + \{H_1,H_2\}_{v,w} {\frac{\partial}{\partial u}} +
\{H_1,H_2\}_{w,u} {\frac{\partial}{\partial v}},
\end{eqnarray}
where the coefficients functions are as given in (\ref{pbh}). Compare the definition of evolutionary vector field $E_{H_1,H_2}$ in (\ref{tNH}) with the one (\ref{cosymEVF}) in the realm of  cosymplectic  framework. The dynamics governing the motion, that is the Nambu-Hamiltonian vector field can be defined as 
\begin{equation} \label{NH-3D}
i_{X_{H_1,H_2}}\mu = dH_1\wedge dH_2 - \nu(H_1)\eta\wedge dH_2 - dH_1\wedge \nu(H_2)\eta, \qquad i_{X_{H_1,H_2}}\eta = 0.
\end{equation}
The time-dependent Nambu-Hamiltonian vector field $X_{H_1,H_2}$ locally looks like the same as the time-independent one in (\ref{HbhVF}). But this time, since the Hamiltonian functions may involve the time variable $t$, in the equations of motion, one may observe some terms depending on $t$ explicitly. If either of the Hamiltonian functions does not depend on the time variable $t$, then this definition reduces to the one in (\ref{NHcf}). Note that, as in the case of the cosymplectic theory, the evolutionary vector field $E_{H_1,H_2}$ can be defined by a summation of the Nambu-Hamiltonian vector field $X_{H_1,H_2}$ and the vertical vector field $\nu$.  

Let us depict now how one can write a Nambu-Hamiltonian system in terms of the differential forms. To this end, recall the associated volume form $\mu$ in (\ref{volume}) and consider the following three form
\begin{equation} \label{mu-H}
\mu_H = \mu - dH_1 \wedge dH_2 \wedge dt,
\end{equation}
on the extended phase space $\mathcal{P}\times \mathbb{R}$. 
The role of this two-form is similar to the role of the Poincar\'{e}--Cartan one-form (\ref{PC}). 
Using $\mu_H$, the Nambu-Hamiltonian vector field can be defined as 
\begin{equation}
i_{X_{H_1,H_2}}\mu_H = 0, \qquad i_{X_{H_1,H_2}}\eta = 0.
\end{equation}
See, for example, \cite{Fe92}. The three-form $\mu_H$ is decomposable. To see this, define the following one-forms
\begin{eqnarray}  \label{beta-3}
\beta^{(1)}:&=&du-\{H_1,H_2\}_{v,w}dt, \quad  \beta^{(2)}:=dv-\{H_1,H_2\}_{u,w}dt, \notag \\
\beta^{(3)}:&=&dw-\{H_1,H_2\}_{u,w}dt.
\end{eqnarray}
It is immediate to check that $\mu_H$ is the exterior product of these one-forms, that is \begin{equation} \label{mu-H2}
\mu_H = \beta^{(1)} \wedge \beta^{(2)} \wedge \beta^{(3)}.
\end{equation} 
In the following subsection, under the light of the present discussion, we shall propose a generalization of the method of Jacobi's last multiplier for the Nambu-Hamiltonian systems.

\subsection{The method of Jacobi's last multiplier for Nambu systems} \label{3jlm}

Now consider the following system of differential equations
\begin{equation}  \label{3D}
\begin{cases}
\dot{x}=f(x,y,z,t),\\
\dot{y}=g(x,y,z,t),  \\ \dot{z}=h(x,y,z,t)
\end{cases}
\end{equation}
depending on three space variables $(x,y,z)$ and a time variable $t$. We introduce the following one-forms
\begin{eqnarray*}
\alpha^{(1)}:&=dx-f(x,y,z,t)dt,\qquad \alpha^{(2)}:=dy-g(x,y,z,t)dt, \\ &\alpha^{(3)}:=dz-h(x,y,z,t)dt.
\end{eqnarray*}
Motivating from the case of $2$-dimensional version presented in the subsection (\ref{JLM2D}), we represent the system (\ref{3D}) by the following three-form \begin{equation} \label{alp-3}
\alpha^{(1)}\wedge \alpha^{(2)}\wedge \alpha^{(3)}= dx\wedge dy \wedge dz - (f dy\wedge dz + g dz \wedge dx + h dx \wedge dy )\wedge dt.
\end{equation}

If the three dimensional system (\ref{3D}) can be written as a Nambu-Hamiltonian form then three-form in (\ref{alp-3}) and the three-form $\mu_H$ in (\ref{mu-H2}) must be multiple of each other at most by a conformal parameter $M$, that is
\begin{equation}\label{M3D}
 M (dx\wedge dy \wedge dz - (f dy\wedge dz + g dz \wedge dx + h dx \wedge dy )\wedge dt)=\mu - dH_1 \wedge dH_2 \wedge dt.
\end{equation}
Let us recall how we have presented the method of Jacobi's last multiplier for the case of $2D$ systems in Subsection(\ref{JLM2D}), and try to mimic all these steps in the present case of Nambu-Hamiltonian systems. The three-form on the right hand side of (\ref{M3D}) is closed, so the one on the left must be closed as well. This gives that 
\begin{equation}
\frac{\partial M}{\partial t}+ f\frac{\partial M}{\partial x}
+g\frac{\partial M}{\partial y}+h\frac{\partial M}{\partial z}
+M\left(\frac{\partial f}{\partial x}+\frac{\partial g}{\partial y}+\frac{\partial h}{\partial z}\right)=0.
\end{equation}
We rewrite this equation as
\begin{equation}\label{JLM3D}
  {\frac{dM}{dt}} +M \hbox{div}X= 0,
\end{equation}
where $X$, in this case, is the vector field generating the dynamics in (\ref{3D}). If the system $\dot{x}=X$ is divergence free with respect to the volume form $\mu$, then it is Nambu-Hamiltonian (\ref{NHcf}) by choosing $M$ as a constant function. If not, the Jacobi's last multiplier $M$ must be non-constant. That is, if the generating vector field $X$ is not divergence free then $M$ can not be a constant function. Assume particularly that, $M$ depends only on the space variables, that is $M=M(x,y,z)$. In the coordinates $(x,y,z)$, the Nambu-Poisson bracket turns out to be 
\begin{equation}
\{F_1,F_2,F_3\}\mu =\frac{1}{M} d F_1 \wedge \nabla F_2 \wedge \nabla F_3.
\end{equation}
This result agrees with the general solution of the Jacobi identity (\ref{jcbv}) and the most general form of the Nambu-Hamilton's equations (\ref{bi-Ham}). In this case, the Jacobi's last multiplier $M$ is a real valued function which makes $MX$ a divergence free vector field. See that, $M$ satisfies 
\begin{equation}  \label{M-3D}
M (f dy\wedge dz + g dz \wedge dx + h dx \wedge dy )= dH_1 \wedge dH_2+\mbox{terms proportional to }\; dt.
\end{equation}

A generalization of the present framework is also possible for the case of time-dependent multiplier $M$. This is achieved by the introduction of three auxiliary functions $\psi$, $\phi$, and $\varphi$  (possibly) depending on the time and the space variables. They are defined through the following equality
\begin{eqnarray}\label{Ham3D}
  M((f-\psi)dy\wedge dz+(g-\phi)dz\wedge dx+(h-\varphi) dx\wedge dy \notag \\= dH_1 \wedge dH_2  + ~ \mbox{terms involving }\; dt.
\end{eqnarray}
It is a matter of a direct calculation to observe that the standard coordinate system $(u,v,w)$ is related with this system as follows
\begin{equation}\label{stcoord}
  M(dx-\psi dt) \wedge (dy - \phi dt) \wedge (dz - \varphi dt)=du \wedge dv \wedge dw.
\end{equation}
This enables us to find a (possibly time-dependent) Poisson transformation determining the standard coordinates of the system. 

\subsection{Three dimensional conformal Nambu-Hamiltonian systems} 

Let $(P,\mu)$ be a three dimensional volume manifold and consider the canonical Nambu-Poisson bracket (\ref{NP-Vol}) associated with $\mu$. Recall that, the flows of a Nambu-Hamiltonian vector field $X_{H_1,H_2}$ preserves the volume form, that is \eqref{NH-vf-mu}. Motivating by the conformal Hamiltonian formalism presented in (\ref{chs}), by loosing the equations \eqref{NH-vf-mu}, we define a conformal Nambu-Hamiltonian vector field $\Gamma^a$ satisfying the identity 
\begin{equation} \label{cNH}
\mathcal{L}_{\Gamma^a}\mu =a\mu
\end{equation}
for some conformal factor $a$. This definition may lead interesting studies both in the geometrical and the mechanical perspectives. 

Scalar multiple $a\mu$ of the volume form $\mu$ is closed. Assuming that it is locally exact, in the standard coordinates $(u,v,w)$, we can write it as the exterior derivative of a two-form 
\begin{equation} \label{tf}
\zeta^{(a_1,a_2,a_3)}=a_1 udv\wedge dw+a_2 v dw\wedge du+a_3wdu\wedge dv.
\end{equation}
Here, we have that the scalar $a$ is $a_1+a_2+a_3$, and $a\mu$ equals to $d\zeta$. 
In this case, we can write (\ref{cNH}) as
\begin{equation}
d(i_{\Gamma^{(a_1,a_2,a_3)}}\mu-\zeta^{(a_1,a_2,a_3)})=0.
\end{equation}
This equation defines $\Gamma^{(a_1,a_2,a_3)}$ modula any exact two-form $dH_1\wedge dH_2$. This enables to write the following identity 
\begin{equation} \label{cNH-Eq}
i_{X_{H_1,H_2} ^{(a_1,a_2,a_3)}}\mu=dH_1\wedge dH_2 +\zeta^{(a_1,a_2,a_3)}. 
\end{equation}
from which we obtain the local formulation of a conformal Nambu-Hamiltonian vector field as follows 
\begin{eqnarray} \label{HbhVF-}
X_{H_1,H_2} ^{(a_1,a_2,a_3)}(u,v,w)&=&  \left(\{H_1,H_2\}_{v,w} +a_1 u\right) {\frac{\partial}{\partial u}} + 
(\{H_1,H_2\}_{w,u} + a_2v){\frac{\partial}{\partial v}} \notag \\&&+ (\{H_1,H_2\}_{u,v} + a_3 w) {\frac{%
\partial}{\partial w}}.
\end{eqnarray}
Note that, we are encoding the conformal parameters $(a_1,a_2,a_3)$ in the notation of the vector field $X_{H_1,H_2} ^{(a_1,a_2,a_3)}$. Let us define a vector field $Z^{(a_1,a_2,a_3)}$ by the equality $$i_{Z^{(a_1,a_2,a_3)}}\mu = \zeta^{(a_1,a_2,a_3)}.$$ Then we see that $X_{H_1,H_2} ^{(a_1,a_2,a_3)}$ is the sum of the Nambu-Hamiltonian vector field $X_{H_1,H_2}$ and $Z$.

\subsection{Example: L\"{u} system} \label{Ex1}

Chaotic L\"{u} system consists of three autonomous first order differential
equations
\begin{equation}  \label{Lu}
\begin{cases}
\dot{x} & =\alpha (y-x), \\
\dot{y} & =\gamma y-xz, \\
\dot{z} & =xy-\beta z,%
\end{cases}%
\end{equation}%
where $\alpha ,\beta $ and $\gamma $ are real constant parameters, and an
overdot represents the derivative with respect to time variable $t$ \cite%
{LuChen,LuChen1}. For $\alpha=36,\;\beta=3$ and $%
\gamma=20$, L\"{u} system is in a chaotic state.

In order to arrive the Hamiltonian formulation of the system, we introduce three auxiliary functions $\psi=-\alpha x$ , $\phi=\gamma y$, and $\varphi=-\beta z $, whereas we take the Jacobi's last multiplier
$M=e^{(\alpha+\beta-\gamma) t} $. The system is in the Nambu-Hamiltonian and the bi-Hamiltonian formulation, see also \cite{EsGhGu16}. To exhibit this, we introduce the Hamiltonian pair
\begin{equation}  \label{Lu5}
H_1(x,y,z) =({\frac{1}{2}}x^2 - \alpha z), \qquad H_2(x,y,z)={\frac{1}{2}}(y^2+z^2) e^{\alpha+\beta-\gamma t}.
\end{equation}
It is now immediate to check that (\ref{Ham3D}) is satisfied. Using (\ref{stcoord}), we compute standard coordinates $(u,v,w)$ as follows
\begin{equation}
  u=xe^{\alpha t}, \qquad v=ye^{-\gamma t}, \qquad w=ze^{\beta t}.
\end{equation}
In this coordinate frame, the L\"{u} system (\ref{Lu}) turns out to be a non-autonomous system of equations
\begin{equation}  \label{Lu-}
\begin{cases}
\dot{u} & =\alpha v e^{(\alpha+\gamma)t}, \\
\dot{v} & =- uw e^{-(\alpha+\beta + \gamma)t}, \\
\dot{w} & =uv e^{(-\alpha+\beta + \gamma)t},
\end{cases}
\end{equation}
whereas the Hamiltonian pair $(H_1,H_2)$ in (\ref{Lu5}) becomes  
\begin{eqnarray} \label{HP-Lu-1}
H_1(u,v,w,t) &=&{\frac{1}{2}}u^2e^{-2\alpha t} - \alpha w e^{-\beta t}\notag  ,\\ H_2(u,v,w,t)&=&
\frac{1}{2}v^2e^{(\alpha+\beta+\gamma) t}+\frac{1}{2}w^2e^{(\alpha-\beta-\gamma) t}.
\end{eqnarray}
The Nambu-Hamiltonian formulation of this system is possible using the geometric framework presented in (\ref{Na3D}). This time, according to the equation 
\begin{equation}
i_{X_{H_1,H_2}}\mu = dH_1\wedge dH_2 - \frac{\partial H_1}{\partial t}dt\wedge dH_2 - dH_1\wedge \frac{\partial H_2}{\partial t}dt, \qquad i_{X_{H_1,H_2}}dt = 0,
\end{equation}
we determine the Nambu-Hamiltonian vector field
$$X_{H_1,H_2}=\alpha v e^{(\alpha+\gamma)t}\frac{\partial}{\partial u}- uw e^{-(\alpha+\beta + \gamma)t}\frac{\partial}{\partial v}+ uv e^{(-\alpha+\beta + \gamma)t}\frac{\partial}{\partial w}. $$ 

An alternative way to study L\"{u} system (\ref{Lu}) is to show that it is possible to recast it in the conformal Nambu-Hamiltonian form (\ref{cNH-Eq}). In this case, the set of conformal parameters are $(-\alpha,\gamma,-\beta)$. We adopt the two-form in (\ref{tf}) to the present case as follows 
\begin{equation}
\zeta^{(-\alpha,\gamma,-\beta)}=-\alpha xdy\wedge dz+ \gamma y dz\wedge dx -\beta zdx\wedge dy.
\end{equation}
By taking $F_1=H_1$ in (\ref{HP-Lu-1}), and by removing the time variable in the definition of $H_2$, we introduce the following time-independent Hamiltonian pair 
\begin{equation}
F_1(x,y,z) =({\frac{1}{2}}x^2 - \alpha z), \qquad F_2(x,y,z)={\frac{1}{2}}(y^2+z^2)
\end{equation}
and obeying the conformal Nambu-Hamiltonian vector field definition in (\ref{HbhVF-}), we compute conformal Nambu-Hamiltonian vector field
\begin{equation}
X_{F_1,F_2} ^{(-\alpha,\gamma,-\beta)}= \alpha (y-x) \frac{\partial}{\partial x} +\gamma y-xz \frac{\partial}{\partial y}+(xy-\beta z)\frac{\partial}{\partial z}. 
\end{equation}
Note that, this vector field exactly generates the L\"{u} system \eqref{Lu}. It is also evident that, this formalization needs either no auxiliary variables or no Jacobi's last multiplier, and the exhibition of both of the Hamiltonian functions and the vector field are time-independent.

\subsection{Example: Qi system} \label{Ex2}
Consider the following reduced for the chaotic Qi system
\begin{equation}  \label{LXsys1}
\begin{cases}
\dot{x} & =y-x+yz, \\
\dot{y} & =\gamma x-xz-y, \\
\dot{z} & =xy-\beta z%
\end{cases}%
\end{equation}
involving cross product terms in each of its arguments \cite{Qi}. We introduce the auxiliary functions $\psi=- x$, $\phi=-y$, $\varphi=-\beta z$ and the Jacobi's last multiplier $M=e^{\beta t}$. Hamiltonian function pair for the system can be choosen as
\begin{equation}
H_{1}=(\gamma x^{2}-y^{2}-(\gamma+1)z^{2})e^{\beta t},\qquad H_{2}=\frac{1}{4(\gamma+1)}(x^{2}+y^{2})-\frac{1}{2}z  \label{LXsys}
\end{equation}
satisfying the defining equation (\ref{Ham3D}). According to the equation (\ref{stcoord}), we compute the standard coordinates as follows
\begin{equation}
u=xe^{t},\qquad v=ye^{t},\qquad w=ze^{\beta t}.  \label{LXsys2}
\end{equation}
In this frame, the Qi system (\ref{LXsys1}) turns out to be non-autonomous 
\begin{equation}
\begin{cases}
\dot{u} & =v+vwe^{-\beta t} \\ \dot{v} & =\gamma u-uwe^{-\beta t} \\ \dot{w} & =uve^{(\beta-2)t}
\label{LXsys3}
\end{cases}
\end{equation}
whereas the Hamiltonian functions becomes 
\begin{eqnarray} \label{HP-Qi-2}
H_{1}&=&(\gamma u^{2}-v^{2})e^{(\beta -2) t}-(\gamma+1)w^{2}e^{-\beta t},\notag \\ H_{2}&=&\frac{1}{4(\gamma+1)}(u^{2}+v^{2})e^{-2t}-\frac{1}{2}w e^{-\beta t}.
\end{eqnarray}
Now, it becomes a matter of direct calculation to show that the Qi system (\ref{LXsys3}) satisfies the time-dependent Nambu formalism presented in (\ref{tNH}), that is 
\begin{equation}
  i_{X_{H_1,H_2}}\mu=dH_1\wedge dH_2 - \frac{\partial H_1}{\partial t} dt \wedge dH_2 - dH_1\wedge \frac{\partial H_2}{\partial t} dt , \qquad i_{X_{H_1,H_2}}dt = 1,
\end{equation}
This results with the Hamiltonian vector field
 \begin{equation}
X_{H_1,H_2} = (v+vwe^{-\beta t}) \frac{\partial}{\partial u} +(\gamma u-uwe^{-\beta t}) \frac{\partial}{\partial v}+(uve^{(\beta-2)t}) \frac{\partial}{\partial w},
\end{equation}
where the Hamiltonian pair $(H_1,H_2)$ in (\ref{HP-Qi-2}) is employed.

Let us now, exhibit the conformal Nambu-Hamiltonian character (\ref{cNH-Eq}) of the Qi system  (\ref{LXsys1}). Bu setting the conformal parameters $(-\alpha,\gamma,-\beta)$, the two-form in (\ref{tf}) turns out to be 
\begin{equation}
\zeta ^{(-1,-1,-\beta)}=-xdy\wedge dz - y dz\wedge dx -\beta zdx\wedge dy.
\end{equation}
Consider the following time-independent Hamiltonian pair 
\begin{equation}
F_{1}=(\gamma x^{2}-y^{2}-(\gamma+1)z^{2}),\qquad F_{2}=\frac{1}{4(\gamma+1)}(x^{2}+y^{2})-\frac{1}{2}z. \end{equation}
See that, $F_{1}=H_1$ in (\ref{HP-Qi-2}), and $F_{2}$ differs $H_2$ by a time-dependent factor.  
and obeying the conformal Nambu-Hamiltonian vector field definition in (\ref{HbhVF-}), we compute the following conformal Hamiltonian operator 
\begin{equation}
X_{F_1,F_2} ^{(-1,-1,-\beta)}= (y-x+yz) \frac{\partial}{\partial x} +(\gamma x-xz-y) \frac{\partial}{\partial y}+(xy-\beta z) \frac{\partial}{\partial z}. 
\end{equation}
It is immediate to see that this time-independent vector field generates the Qi system  (\ref{LXsys1}).

\section{Conclusion and Discussions}
In this work, we have focused on the non-canonical Hamiltonian realizations of $2D$ and $3D$ dynamical systems. For $2D$ systems, in subsection (\ref{JLM2D}), time-dependent Hamiltonian formulations, obtained by the method of Jacobi's last multiplier, are studied in the framework of cosymplectic geometry. To deal with affine terms in $2D$, we have employed the theory of conformal Hamiltonian dynamics. A host-parasite model has been extensively studied in the subsection (\ref{hpm}). In the table (\ref{table}), we have presented Hamiltonian formulations of some $2D$ biological systems.   

For $3D$ systems, we have first recalled some basics of Poisson and Nambu-Poisson manifolds. The equation (\ref{NHcf}) is representing the Nambu-Hamiltonian systems in a covariant way. We have generalized the covariant formulation (\ref{NHcf}) of the Nambu-Hamiltonian systems for the non-autonomous cases by the introduction of the equation (\ref{tNH}). Accordingly, in the subsection (\ref{3jlm}), the method of Jacobi's last multiplier has been presented for $3D$ Nambu systems. A conformal generalization of Nambu-Hamiltonian geometry has been introduced in (\ref{cNH-Eq}). Two examples have been provided, namely L\"{u} system and Qi system in subsections (\ref{Ex1}) and (\ref{Ex2}), respectively.

\section*{Funding statement}

The research of Partha Guha is partially supported by FAPESP through Instituto 
de Fisica de S\~ao Carlos, Universidade de Sao Paulo with grant number 2016/06560-6.

\end{document}